\documentclass[notitlepage,11pt]{article}
\usepackage{amssymb,amsfonts,amsmath,geometry
}

\catcode`\@=12
\usepackage[notref,notcite]{showkeys}

\newtheorem{Theorem}{Theorem}
{Lemma}
{Proposition}
{Corollary}
\newtheorem{Remark}
{Remark}

\usepackage{color}

\def\R{\mathbb{R}}

\def\irn{\int\limits_{\R^n}}

\def\div{{\rm div}}

\def\Dshalf{\left(-\Delta\right)^{\!\frac{s}{2}}}  

\def\proof{\noindent{\textbf{Proof. }}}
\def\QED{\hfill {$\square$}\goodbreak \medskip}

\usepackage{mathtools}

\begin{document}

\title{A weighted estimate for generalized harmonic extensions}

\author{Roberta Musina
\footnote{Dipartimento di Matematica ed Informatica, Universit\`a di Udine,
via delle Scienze, 206 -- 33100 Udine, Italy. Email: {roberta.musina@uniud.it} 
}~ 
and
Alexander I. Nazarov
\footnote{
St.Petersburg Department of Steklov Institute, Fontanka 27, St.Petersburg, 191023, Russia, 
and St.Petersburg State University, 
Universitetskii pr. 28, St.Petersburg, 198504, Russia. E-mail: al.il.nazarov@gmail.com
}
}

\date{}

\maketitle

\begin{abstract}
{\small We prove some weighted $L_p$ estimates for 
generalized harmonic extensions in the half-space.}

\vskip0.5cm


\medskip\noindent
\textbf{2010 Mathematics Subject Classfication:}  	35A23;	42B35. 
\end{abstract}

\vskip1cm

Let $u=u(x)$ be a ``good'' function in $\R^n$. 
Denote by ${\mathbb P}u=({\mathbb P}u)(x,y)$ its 
harmonic extension to the half-space $\R^{n+1}_+\equiv \R^n\times(0,\infty)$,
$$
({\mathbb P}u)(x,y)=
\frac{\Gamma(\frac {n+1}2)}{\pi^{\frac {n+1}2}}\,\irn  u(\xi)\cdot\frac{y\,d\xi}{\big(|x-\xi|^2+y^2\big)^{\frac {n+1}2}},\qquad x\in\R^n, y>0.
$$
By elementary convolution estimates, the linear map $\mathbb P:u\mapsto ({\mathbb P}u)(\,\cdot\,,y)$ is non-expansive in $L_p(\R^n)$
for any   $p\in[1,\infty]$, that is, $\|({\mathbb P}u)(\,\cdot\,,y)\|_p\le \|u\|_p$  for any $y>0$.
%

In the breakthrough paper \cite{CaSi}, Caffarelli and Silvestre introduced, for any
 $s\in(0,1)$, the following
generalized {\em $s$-harmonic extension}
$u~\mapsto ~{\mathbb P_{\!s}}u$,
$$
({\mathbb P_{\!s}}u)(x,y)=c_{n,s}
\,\irn  u(\xi)\cdot\frac{y^{2s}\,d\xi}{\big(|x-\xi|^2+y^2\big)^{\frac {n+2s}2}}\,,
\quad c_{n,s}=\frac{\Gamma\big(\frac{n+2s}{2}\big)}{\pi^{\frac{n}{2}}\Gamma(s)}
$$
so that the classical harmonic extension is recovered for $s=\frac 12$.
One of the main results in  \cite{CaSi} states that 
the $L_2$-norm of $\Dshalf u=\mathcal F^{-1}\big[|\xi|^s\mathcal F[u]\big]$ on $\R^n$
(here $\mathcal F$ is the Fourier transorm in $\R^n$) coincides, up to a constant that depends only on $s$, 
with
some weighted $L_2$-norm of $|\nabla({\mathbb P_{\!s}}u)|$ on $\R^{n+1}_+$.

Notice that for arbitrary $y>0$, the kernel
\begin{equation}
\label{eq:s_extension}
\mathcal P_{\!s}(x,y)= {\frac{\Gamma\big(\frac{n+2s}{2}\big)}{\pi^{\frac{n}{2}}\Gamma(s)}}\,
\frac{y^{2s}}{\big(|x|^2+y^2\big)^{\frac {n+2s}2}}
\end{equation}
has unitary $L_1$-norm, thus
the linear map $u\mapsto ({\mathbb P_{\!s}}u)(\cdot,y)$ is  non-expanding in $L_p(\R^n)$
as well. In particular, we have 
$$
\irn |({\mathbb P_{\!s}}u)(\,\cdot\,,y)|^p\,dx \le \irn |u|^p\,dx \quad \text{for any $s\in(0,1)~,~~ y>0~,~~ p\in[1,\infty)$.}
$$
We are interested in similar results for weighted $L_p$-norms.
More precisely, we deal with inequalities of the form
\begin{equation}\label{more}
\irn \frac{|({\mathbb P_{\!s}}u)(x,y)|^p}{\big(|x|^2+y^2\big)^{\alpha}}\,dx\le C_p\irn \frac{|u(x)|^p}{\big(|x|^2+y^2\big)^{\alpha}}\,dx
\end{equation}
where $C_p>0$ does not depend on $y, u$. These inequalities seems to be new even in the classical
case $s=\frac12$.

The next statement is crucially used in \cite{MN-SB}.

\begin{Theorem}
Let $s\in(0,1)$, $\alpha\ge 0$.

\begin{itemize}
\item[$i)$] If $p=1$, The inequality (\ref{more}) 
holds  if and only if $\alpha\le \frac n2+s$.
\item[$ii)$] For arbitrary $1<p<\infty$, the inequality (\ref{more}) 
holds  if and only if $\alpha<\frac n2+sp$.
\end{itemize}
\end{Theorem}

\proof
Take a measurable function $u$, an arbitrary $y>0$, and put $u^y(x)=u(yx)$. By dilation, we have
$({\mathbb P_{\!s}}u)(x,y)=({\mathbb P_{\!s}}u^y)(\frac{x}{y},1)$. Thus it suffices to prove 
(\ref{more}) for $y=1$. 

In case $p=1$, we rewrite the inequality
\begin{equation}\label{one}
\irn \frac{|({\mathbb P_{\!s}}u)(x,1)|}{\big(|x|^2+1\big)^{\alpha}}\,dx\le C_1\irn \frac{|u(x)|}{\big(|x|^2+1\big)^{\alpha}}\,dx
\end{equation}
in the  form
\begin{equation}
\label{eq:equivalent}
\irn \Big|\irn \frac{u(\xi)}{\big(|\xi|^2+1\big)^\alpha}\,
\frac{c_{n,s}}{(|x-\xi|^2+1)^{\frac{n+2s}{2}}}
\frac{\big(|\xi|^2+1\big)^\alpha}{\big(|x|^2+1\big)^\alpha}\,d\xi\Big|dx\le 
C_1\irn\frac{|u(\xi)|}{\big(|\xi|^2+1\big)^\alpha}\,d\xi,
\end{equation}
to make evident that we are indeed estimating the norm of the transform
$$
v\mapsto \mathbb Lv~,\quad 
({\mathbb L}v)(x)=\irn v(\xi)\,
\frac{c_{n,s}}{(|x-\xi|^2+1)^{\frac{n+2s}{2}}}
\frac{\big(|\xi|^2+1\big)^\alpha}{\big(|x|^2+1\big)^\alpha}\,d\xi
$$
as a linear operator $L_1(\R^n)\to L_1(\R^n)$. We use the duality $L_1(\R^n)'=L_\infty(\R^n)$,
that gives
\begin{equation}\label{supremum}
\|{\mathbb L}\|_{L_1\to L_1}=\sup_{\xi\in\R^n}\,
\irn \frac {c_{n,s}}{(|x-\xi|^2+1)^{\frac {n+2s}2}}\,\frac {\big(|\xi|^2+1\big)^{\alpha}}{(|x|^2+1)^{\alpha}}\,dx.
\end{equation}
If $\alpha>\frac n2 +s$, then the supremum in (\ref{supremum}) is evidently infinite. If $\alpha\le\frac n2 +s$ then easily
$$
 \int\limits_{|x|\ge|\xi|/2} \frac {c_{n,s}}{(|x-\xi|^2+1)^{\frac {n+2s}2}}\,
 \frac {\big(|\xi|^2+1\big)^{\alpha}}{\big(|x|^2+1\big)^{\alpha}}\,dx\le \irn 2^{2\alpha}\mathcal P_{\!s}(x-\xi,1)\,dx= 2^{2\alpha}.
$$
Further, $|x|\le|\xi|/2$ implies $|x-\xi|\ge|\xi|/2$ and $ |x-\xi|\ge|x|$. Therefore,
\begin{eqnarray*}
 \int\limits_{|x|\le|\xi|/2} \frac {c_{n,s}}{(|x-\xi|^2+1)^{\frac {n+2s}2}}\,\frac {\big(|\xi|^2+1\big)^{\alpha}}{\big(|x|^2+1\big)^{\alpha}}\,dx
 &\le &
 \int\limits_{|x|\le|\xi|/2} \frac {2^{2\alpha}c_{n,s}\,dx}{(|x-\xi|^2+1)^{\frac {n+2s}2-\alpha}\big(|x|^2+1\big)^{\alpha}}\\
 &\le &
 \irn 
  2^{2\alpha}\mathcal P_{\!s}(x,1)\,dx= 2^{2\alpha}.
\end{eqnarray*}
We can conclude that $C_1=\|{\mathbb L}\|_{L_1\to L_1}<\infty$, and $i)$ is proved.
\medskip

Next, we take $p>1$. To handle the case  $\alpha\ge \frac{n}{2}+sp$ we notice that the function 
$$
\overline{u}(x):=\frac {(|x|^2+1)^{\frac {2\alpha-n}{2p}}}{\log(|x|^2+2)}
$$
satisfies
$$
\irn \frac{|\overline{u}(x)|^p}{\big(|x|^2+1\big)^{\alpha}}\,dx=\irn \frac{dx}{(|x|^2+1)^{\frac n2}{\log^p(|x|^2+2)}}<\infty,
$$
On the other hand,  for any arbitrary $x\in\R^n$ we have
$$
\irn \mathcal P_{\!s}(x-\xi,1)\overline{u}(\xi)\,d\xi>\irn \frac {C(x)\,d\xi}{(|\xi|^2+1)^{\frac n2}\log(|\xi|^2+2)}\,,
$$
and the last integral diverges. Thus, for $p>1$ and $\alpha\ge\frac n2+sp$ the inequality (\ref{more}) does not hold with a finite constant $C$
in the right hand side.

If $\alpha<\frac n2 +sp$, we use H\"older's inequality to get
$$
|({\mathbb P_{\!s}}u)(x,1)|\le \Big(\irn \mathcal P_{\!s}(x-\xi,1)\,\frac {|u(\xi)|^p}{\big(|\xi|^2+1\big)^{\beta}}\,d\xi\Big)^{\frac 1p}
\Big(\irn \mathcal P_{\!s}(x-\xi,1)\big(|\xi|^2+1\big)^{\frac {\beta}{p-1}}\,d\xi\Big)^{\frac {p-1}p},
$$
where $\beta:=\max\{\alpha-\frac n2-s,0\}<s(p-1)$.

If $\alpha\le\frac n2 +s$ then $\beta=0$ and  the last integral equals $1$. In this case we obtain
\begin{equation}\label{le1}
\irn \frac{|({\mathbb P_{\!s}}u)(x,1)|^p}{\big(|x|^2+1\big)^{\alpha}}\,dx\le 
\irn {\mathbb L}
\bigg[\frac {|u(\cdot)|^p}{\big(|\cdot|^2+1\big)^{\alpha}}\bigg](x)\,dx
\le \|\mathbb L\|_{L_1\to L_1}
\irn 
\frac {|u(x)|^p}{\big(|x|^2+1\big)^{\alpha}}\,dx\,,
\end{equation}
and (\ref{more}) follows from the first part of the proof. 

If $\frac{n}{2}+s<\alpha<\frac{n}{2}+sp$, we estimate
$$
\begin{aligned}
\irn \frac{|({\mathbb P_{\!s}}u)(x,1)|^p}{\big(|x|^2+1\big)^{\alpha}}\,dx \le 
\Big(\irn \irn \mathcal P_{\!s}(x-\xi,1)&~\,\frac {\big(|\xi|^2+1\big)^{\frac {n+2s}2}}{\big(|x|^2+1\big)^{\frac {n+2s}2}}\,\frac {|u(\xi)|^p}{\big(|\xi|^2+1\big)^{\alpha}}\,d\xi dx\Big) \\
\times & \Big(\sup_{x\in\R^n}\irn \frac {c_{n,s}}{(|x-\xi|^2+1)^{\frac {n+2s}2}}\,\frac {\big(|\xi|^2+1\big)^{\frac {\beta}{p-1}}}{\big(|x|^2+1\big)^{\frac {\beta}{p-1}}}\,d\xi\Big)^{p-1}.
\end{aligned}
$$
If we prove that the last supremum is finite then (\ref{more}) again follows from the first statement of the present theorem. We have
$$
\int\limits_{|\xi|\le2|x|} \frac {c_{n,s}}{(|x-\xi|^2+1)^{\frac {n+2s}2}}\,\frac {\big(|\xi|^2+1\big)^{\frac {\beta}{p-1}}}{\big(|x|^2+1\big)^{\frac {\beta}{p-1}}}\,d\xi\le 
2^{\frac {2\beta}{p-1}}\irn \mathcal P_{\!s}(x-\xi,1)\,d\xi= 2^{\frac {2\beta}{p-1}}.
$$
Further, $|\xi|\ge 2|x|$ implies $|x-\xi|\ge|\xi|/2$. Therefore, from $\beta<s(p-1)$ we get
$$
\int\limits_{|\xi|\ge 2|x|} \frac {c_{n,s}}{(|x-\xi|^2+1)^{\frac {n+2s}2}}\,\frac {\big(|\xi|^2+1\big)^{\frac {\beta}{p-1}}}{\big(|x|^2+1\big)^{\frac {\beta}{p-1}}}\,d\xi
\le
 \int\limits_{|\xi|\ge 2|x|} \frac {2^{\frac {2\beta}{p-1}}c_{n,s}\,d\xi}{(|x-\xi|^2+1)^{\frac n2 +s-\frac {\beta}{p-1}}} \le C(n,s,\beta,p),
$$
and the proof of (\ref{more}) is complete.
\QED

The following statement partially solves the problem whether the map $u\mapsto ({\mathbb P_{\!s}}u)(\cdot,y)$ is non-expanding in weighted $L_p(\R^n)$.

\begin{Theorem}
Let $s\in(0,1)$.

\begin{itemize}
\item[$i)$] If $0\le\alpha\le \frac n2-s$ then for arbitrary $1\le p<\infty$ the best constant in (\ref{more}) is $C_p=1$.
\item[$ii)$] If $\alpha>\frac n2$ then the best constant $C_p$  in (\ref{more}) is greater than $1$,
at least for $p$ close to $1^+$.
\end{itemize}
\end{Theorem}

\begin{Remark}
 We conjecture that the statement $ii)$ holds for all $1\le p<\infty$. The value of $C_p$ for $\frac n2-s<\alpha\le \frac n2$ is a completely open problem.
\end{Remark}

\proof
We again suppose $y=1$. 
\medskip

Firstly, we prove $i)$  in case $p=1$.
It has been proved in \cite{CaSi} that the function 
\begin{equation}\label{CS_ext}
\omega(\xi,y)=\irn \mathcal P_{\!s}(\xi-x,y)\,\frac {dx}{\big(|x|^2+1\big)^{\alpha}}
\end{equation}
solves the following boundary value problem in ${\mathbb R}^{n+1}_+$,
\begin{equation}\label{equation}
 -\div(y^{1-2s}\nabla \omega)=0; \qquad \omega(\xi,0)=\big(|\xi|^2+1\big)^{-\alpha}.
\end{equation}
Consider the barrier function $\widetilde{\omega}(\xi,y)=\big(|\xi|^2+y^2+1\big)^{-\alpha}$. A direct computation gives
$$
-\div(y^{1-2s}\nabla \widetilde{\omega})=2\alpha y^{1-2s}\widetilde{\omega}^{1+\frac 2{\alpha}}\big((n-2s+2)+(n-2s-2\alpha)(|\xi|^2+y^2)\big)\ge 0
$$
because of the assumption on $\alpha$.
Since $\widetilde\omega(\xi,0)=\omega(\xi,0)$, we have that
$\omega\le\widetilde{\omega}$ in $\R^{n+1}_+$ by the maximum principle. In particular,
$$
(|\xi|^2+1)^\alpha\omega(\xi,1)<\Big(\frac{|\xi|^2+1}{|\xi|^2+2}\Big)^\alpha<1.
$$ 
Therefore, the supremum in (\ref{supremum}) does not exceed $1$, and thus 
the best constant in (\ref{one}) is $C_1=\|{\mathbb L}\|_{L_1\to L_1}\le 1$.

Since $\|{\mathbb L}\|_{L_1\to L_1}\le 1$, the inequalities in (\ref{le1})
readily give $C_p\le 1$, for any $p\ge 1$. 

Finally, to prove that $C_p=1$ if $\alpha\le \frac{n}{2}-s$, it suffices to consider the sequence 
$u(\varepsilon x)$, where $u\in{\cal C}^\infty_0({\mathbb R}^n)$, $u\ge0$, is a fixed nontrivial function, and
then to push $\varepsilon$ to $0$. The proof of $i)$ is complete.

\medskip

To prove $ii)$ consider the function $v(x)=\big(|x|^2+1\big)^{-\alpha}$. Clearly $v\in L_1({\mathbb R}^n)$ and
$$
\irn ({\mathbb P}v)(x)\,dx=\irn \mathcal P_{\!s}(\xi-x,1)\,dx\irn v(\xi)\,d\xi=\irn v(\xi)\,d\xi.
$$
Since  
$$
({\mathbb P}v)(0)=\irn \mathcal P_{\!s}(\xi,1)v(\xi)\,d\xi<\max v(\xi)=v(0),
$$
there exists a point $\xi$ such that $({\mathbb P}v)(\xi)>v(\xi)$. Therefore, the supremum in (\ref{supremum}) is greater then $1$, and the best constant in (\ref{one}) is 
$C_1=\|{\mathbb L}\|_{L_1\to L_1}> 1$. By continuity,  the best constant in (\ref{more}) is greater than $1$
for $p$ sufficiently close to $1$.
\QED

\bigskip
{\bf Acknowledgements}. 
The first author is partially supported by Miur-PRIN project 2015KB9WPT\_001 and PRID project {\em VAPROGE}. The second author is partially supported by  RFBR grant 17-01-00678. 
We are grateful to N. Filonov for the hint to the proof of the statement $ii)$ of Theorem 2.

\end{document}